\documentclass[a4paper,11pt]{article}

\scrollmode
\usepackage{amsmath,amssymb,amsfonts,graphicx,color,fancyhdr}
\usepackage[latin1]{inputenc}
\usepackage[pdftex, colorlinks=true, pdfstartview=FitV, linkcolor=blue, citecolor=blue, urlcolor=blue]{hyperref}

\newtheorem{thm}{Theorem}[section]
\newtheorem{prop}[thm]{Proposition}
\newtheorem{lem}[thm]{Lemma}
\newtheorem{df}[thm]{Definition}

\newtheorem{rem}[thm]{Remark}
\newtheorem{cor}[thm]{Corollary}

\def\be#1 {\begin{equation} \label{#1}}
\newcommand{\ee}{\end{equation}}
\def\dem {\noindent {\bf Proof : }}

\def\sqw{\hbox{\rlap{\leavevmode\raise.3ex\hbox{$\sqcap$}}$%
\sqcup$}}
\def\findem{\ifmmode\sqw\else{\ifhmode\unskip\fi\nobreak\hfil
\penalty50\hskip1em\null\nobreak\hfil\sqw
\parfillskip=0pt\finalhyphendemerits=0\endgraf}\fi}

\newcommand{\mb}{\medskip\noindent}
\newcommand{\gb}{\bigskip\noindent}
\newcommand{\R}{\mathbb R}

\newcommand{\N}{\mathbb N}

\newcommand{\C}{\mathbb C}

\def\Xint#1{\mathchoice
   {\XXint\displaystyle\textstyle{#1}}%
   {\XXint\textstyle\scriptstyle{#1}}%
   {\XXint\scriptstyle\scriptscriptstyle{#1}}%
   {\XXint\scriptscriptstyle\scriptscriptstyle{#1}}%
   \!\int}
\def\XXint#1#2#3{{\setbox0=\hbox{$#1{#2#3}{\int}$}
     \vcenter{\hbox{$#2#3$}}\kern-.5\wd0}}

\def\aver#1{\Xint-_{#1}}

\newcommand{\M}{{\mathcal M}}
\title{Abstract framework for John Nirenberg inequalities and applications to Hardy spaces}
\author{Fr\'ed\'eric Bernicot\\ CNRS - Laboratoire Paul Painlev\'e \\ Universit\'e
Lille 1\\59655 Villeneuve d'Ascq Cedex, France\\frederic.bernicot@math.univ-lille1.fr \and Jiman Zhao \footnote{supported by SRF for ROCS, SEM, China, and by NNSF of China No. 10871048 and No. 10931001.} \\ School of Mathematical Sciences, Beijing Normal University \\
Key Laboratory of Mathematics and Complex Systems, Ministry of
Education, \\ Beijing 100875, P.R. China \\ jzhao@bnu.edu.cn}
\date{March 1, 2010}

\begin{document}

\maketitle

\begin{abstract} In this paper, we develop an abstract framework for John-Nirenberg inequalities associated to BMO-type spaces. This work can be seen as the sequel of \cite{BZ}, where the authors introduced a very general framework for atomic and molecular Hardy spaces. Moreover, we show that our assumptions allow us to recover some already known John-Nirenberg inequalities. We give applications to the atomic Hardy spaces too.
\end{abstract}

\mb {\bf Key-words:} BMO spaces ; John-Nirenberg inequalities ; Hardy spaces.

\mb {\bf MSC:} 42B20 ; 46E30.

\tableofcontents

\section{Introduction}

\gb The first BMO space (space of functions satisfying a {\it Bounded Mean Oscillation}) was originally introduced by F. John and L. Nirenberg in \cite{JN}. This space naturally arises as the class of functions whose deviation from their means over cubes is bounded. This space is strictly including the $L^\infty$ space and is a good extension of the Lebesgue spaces scale $(L^p)_{1<p<\infty}$ for $p\to \infty$ from a point of view of Harmonic Analysis. For example, it plays an important role in boundedness of Calder\'on-Zygmund operators, real interpolation, Carleson measure, study of paraproducts,  ... Moreover the BMO space can be characterized by duality as the dual space of the Coifman Weiss space $H^1$. This observation was annouced by C. Fefferman in \cite{F} and then proved in \cite{FS}.

\gb Here we are interested in one of the most important properties
of the BMO space : the so-called John-Nirenberg inequality (see
\cite{JN}). This deep property describes the exponential
integrability of the oscillations for a BMO-function. More
precisely, for $Q$ a ball of the Euclidean space $\R^n$ then a
function $f\in BMO$ satisfies
$$ \left|\left\{x\in Q,\ \left|f-\aver{Q} f \right| >\lambda\right\}\right| \leq c_1|Q|e^{-c_2\lambda/\|f\|_{BMO}}$$
for some constants $c_1,c_2$ only dependent on the dimension $n$. Consequently the oscillation $f-\aver{Q} f $, which was initially supposed to belong to $L^1(Q)$ (by the definition of BMO), is indeed exponentially integrable on $Q$. The BMO-norm arise to a self-improvement of the integrability of the oscillation. \\
The first consequence of such inequalities was the equivalence bewteen the spaces $BMO_{q}$ for $q\in (1,\infty)$ (where $BMO_q$ is based on a control of the oscillations in $L^q$ norm). A second important consequence concerns the Hardy space $H^1$. Since the duality results $(H^1)^*=BMO$, it comes that the Hardy space defined by atomic decomposition with $p$-atoms, does not depend on the exponent $p\in (1,\infty)$ (see the work of R. Coifman and G. Weiss about Hardy spaces \cite{CW}). \\
We just detail these properties as our aim will be to extend these ones in an abstract framework. However BMO space has many properties and the use of this space in real Harmonic Analysis has given rise to many works (related to Calder\'on-Zygmund operators, Carleson measures, $T(1)$-theorems, ...).

\gb Unfortunately, there are situations where the John-Nirenberg space $BMO$ is not the
 right substitute to $L^\infty$ and there have been recently numerous works
 whose the goal is to define an adapted BMO space according to the context (see \cite{DY,HM} ...).
 For example the classical space $BMO$ is not well adapted to operators such as the Riesz transform
  on Riemannian manifolds. That is why in \cite{HM}, S. Hofmann and S. Mayboroda develop theory of Hardy and BMO spaces associated to a second order divergence form elliptic operators, which also including the corresponding John-Nirenberg inequality. In the recent works \cite{DY} and \cite{DY1}, X. T. Duong and L. Yan studied some new BMO type spaces and proved an associated version of John-Nirenberg inequality on these spaces (with duality results).

\gb In \cite{BZ}, the authors have developped an abstract framework for atomic Hardy spaces (and proved some results about interpolation with Lebesgue spaces). Without more precised assumptions, it seems to be impossible to get a precise characterization of the dual space of them by a BMO-type space ; although these last ones are well-defined. This framework permit to cover the classical space $BMO$ and those defined in \cite{DY} and \cite{HM}.

\gb The aim of this article is to continue the study in this abstract framework of BMO-type spaces and to describe assumptions implying John-Nirenberg inequalities associated to these new spaces. We refer to \cite{JM,J,BM} for recent works on general self-improvement properties of some inequalities. Here we are specially interested in John-Nirenberg inequalities. That is why our point of view is based on the associated Hardy spaces and we are looking to Assumptions related to this approach.\\
In detail, our paper is organized as follows: \\
In Section \ref{sec:def}, we define our framework of Hardy and BMO spaces and then state our main results concerning John-Nirenberg inequalities (see Theorems \ref{thm} and \ref{thm2}). We postpone theirs proofs to Subsection \ref{subsec:proof}. In Section \ref{sec:check}, we check that our assumptions are reasonable and that our results
generalize some already known particular cases such as the John-Nirenberg
inequalities in \cite{JN}, \cite{FS}, and \cite{DY}. In section \ref{sec:appli}, we present an application of the John-Nirenberg inequalities to the Hardy spaces.

\section{An abstract framework for John-Nirenberg inequalities} \label{sec:def}

\subsection{Hardy and BMO spaces} \label{subsec:def}

Let $(X,d,\mu)$ be a space of homogeneous type. There exist constants $A,\delta>0$ such that for all $x\in X$, $r>0$ and $t\geq 1$
\be{homogene} \mu(B(x,tr)) \leq A t^{\delta}\mu(B(x,r)), \ee where
$B(x,r)$ is the open ball with center $x\in X$ and radius $r>0$. We
call $\delta$ the homogeneous dimension of $X$. For $Q$ a ball, and
$i\geq 0$, we write $C_i(Q)$ the scaled corona around the ball $Q$~:
$$ C_i(Q):=\left\{ x,\ 2^{i} \leq 1+\frac{d(x,c(Q))}{r_Q} < 2^{i+1} \right\},$$
where $r_Q$  is the radius of the ball $Q$ and $c(Q)$ its center.
 Note that $C_0(Q)$ corresponds to the ball $Q$ and $C_i(Q) \subset 2^{i+1}Q$ for $i\geq 1$, where $\lambda Q$ is as usual the ball with center $c(Q)$ and radius $\lambda r_Q$. For $p\in[1,\infty]$, we denote by $L^p=L^p(X)$ the Lebesgue space.

\mb Let us denote by $\mathcal{Q}$ the collection of all balls~:
$$ \mathcal{Q}:= \left\{ B(x,r),\ x\in X, r>0 \right\}.$$ We write ${\mathcal M}$ for the Hardy-Littlewood maximal operator and for $p\in [1,\infty)$, we denote its $L^p$-version by
$$ {\mathcal M}_p(f)(x) := {\mathcal M}(|f|^p)(x) ^{1/p}.$$
Let $\mathbb{B}:=(B_Q)_{Q\in \mathcal{Q}}$ be a collection of
$L^2$-bounded linear operator, indexed by the collection
$\mathcal{Q}$. We write $A_Q=Id-B_Q$ and $B_Q^*$ for its adjoint
operator. We assume that these operators $B_Q$ are uniformly bounded is some Lebesgue space: there exist two exponents $p_1<p_0$ belonging to $(1,\infty]$ and a constant $0<A'<\infty$ so that~: for all $p\in [p_1,p_0]$
\be{operh} \forall f\in L^p ,\ \forall  Q \textrm{ ball}, \qquad
\|B_Q(f) \|_{L^p} \leq A'\|f\|_{L^p}. \ee

\mb In the rest of the paper, we allow the constants to depend on
$A$, $A'$ and $\delta$.

\mb For convenience, we first recall the definition of {\it atoms},
{\it molecules} and the corresponding Hardy spaces introduced in
\cite{BZ}.

\begin{df}(\cite{BZ}) Let $\epsilon>0$ and $p\in[p_1,p_0]$ be fixed parameters. A function $m\in L^{1}_{loc}$ is called a $(\epsilon,p)$-molecule associated to a ball $Q$ if there exists a real function $f_Q$ such that
$$m=B_Q(f_Q),$$
with
$$\forall i\geq 0, \qquad  \|f_Q\|_{L^p(C_i(Q))} \leq \left(\mu(2^{i}Q)\right)^{-1+1/p} 2^{-\epsilon i}.$$
We call $m=B_Q(f_Q)$ a $p$-atom if in addition we have $supp(f_Q)
\subset Q$. So a $p$-atom is exactly an $(\infty,p)$-molecule.
\end{df}

\begin{df} (\cite{BZ}) A measurable function $h$ belongs to the molecular
Hardy space $H^1_{p,\epsilon,mol}$ if there exists a decomposition~:
$$h=\sum_{i\in\N} \lambda_i m_i  \qquad \mu-a.e, $$
where for all $i$, $m_i$ is an $(\epsilon,p)$-molecule and $\lambda_{i}$
are real numbers satisfying
$$\sum_{i\in \N} |\lambda_i| <\infty. $$
%\qquad \textrm{and} \qquad \sum_{i} \left| \lambda_i m_i\right|<\infty \qquad \mu-a.e
We define the norm~:
$$\|h\|_{H^1_{p,\epsilon,mol}}:= \inf_{h=\sum_{i\in\N} \lambda_i m_i} \sum_{i} |\lambda_i|.$$
Similarly we define the atomic space $H^1_{p,ato}$ replacing
$(\epsilon,p)$-molecules by $p$-atoms.
\end{df}

\mb We refer the reader to \cite{BZ} for the first work introducing these notations and these abstract point of view for Hardy spaces (and to \cite{BB} for extension in considering Hardy-Sobolev spaces). We emphasize that in \cite{BZ} and \cite{B}, interpolation results are described. Moreover in \cite{BZ} the authors have described weak results concerning duality results $H^1-BMO$. In the following, we define the BMO spaces.

\begin{df} For $q\in[p_0',p_1']$, a function $f\in L^q$ belongs to the space $Bmo_{q}$
if
$$ \left\|f \right\|_{BMO_q} := \sup_{Q \textrm{ ball}} \left(\frac{1}{\mu(Q)}\int_{Q} \left| B_Q^*(f) \right|^q d\mu \right)^{1/q } <\infty,$$
where we denote $B_Q^*$ for the adjoint operator of $B_Q$. We write $BMO_q$ for the completion of $Bmo_q$ for the corresponding norm.
\end{df}

\mb We will see that it could be interesting to define other ``molecular'' BMO spaces as follows:

\begin{df} For $\epsilon>0$ and $q\in[p_0',p_1']$, a function $f\in L^q$ belongs to the molecular space $Bmo_{\epsilon,q}$ if
$$ \left\|f \right\|_{BMO_{\epsilon,q}} := \sup_{Q \textrm{ ball}}\ \sup_{j\geq 0}\ 2^{j\epsilon}\left(\frac{1}{\mu(2^j Q)}\int_{C_j(Q)} \left| B_Q^*(f) \right|^q d\mu \right)^{1/q } <\infty.$$
\end{df}

\begin{rem} Obviously we have $BMO_{\epsilon,q} \hookrightarrow BMO_{q}$. The question of a reverse property is open in such an abstract framework.
\end{rem}

\mb After these definitions, we can state our main results in the following subsection.

\subsection{John-Nirenberg inequalities}

First we have to make some Assumptions on the operators $B_Q$ in order to get John-Nirenberg inequalities. Let us first describe the required properties.

\gb {\bf Assumptions.} We set $q_0=p_0'$ and $q_1=p_1'$ such that $1\leq q_0<q_1<\infty$.

\mb We assume $L^{q_0}-L^{q_1}$ off-diagonal decay for operators
$A_Q^*$. There exists coefficients $\gamma_j$ such that for all ball
$Q$
\begin{equation} \label{assum}
 \left(\frac{1}{\mu(Q)}\int_{Q} \left|A_Q^*(h) \right|^{q_1} d\mu \right)^{1/q_1} \lesssim \sum_{j} \gamma_j
\left(\frac{1}{\mu(2^{j}Q)}\int_{C_j(Q)} \left|h \right|^{q_0} d\mu
\right)^{1/q_0},
%\tag{${\mathcal O}$}
\end{equation}
with the property
\be{assum2} \sum_j \gamma_j <\infty. \ee

\begin{rem} \label{rem:inter} Let us note that off-diagonal decay (\ref{assum}) is slighty stronger than the following maximal inequality~: for all $x\in X$
$$ \sup_{Q\ni x}  \left(\frac{1}{\mu(Q)}\int_{Q} \left|A_Q^*(h) \right|^{q_1} d\mu \right)^{1/q_1} \lesssim {\mathcal M}_{q_0}(|f|)(x).$$
This comparing between maximal operators is exactly the assumption, required in our previous works \cite{BZ,B}, in order to get interpolation results between our Hardy spaces and $L^{p_0}$. It is interesting that almost the same assumption seems to be important for both interpolation results and John-Nirenberg inequalities. \\
Specially, we move the reader to Corollary 5.8 of \cite{BZ} for obtaining adapted version of {\it Fefferman-Stein inequalities}.
\end{rem}

%\mb We suppose that the operators $(B_Q)_Q$ commutte : $B_QB_R=B_RB_Q$ for all balls $Q$and $R$ (which is equivalent to the commutativity of operators $A_Q$).

\mb Moreover we make the following assumption~: there exists a
constant $c$ such that
\begin{equation} \label{ass:cont}
\sup_{\genfrac{}{}{0pt}{}{R \textrm{ ball}}{R\subset 2Q}}
\left(\frac{1}{\mu(R)}\int_{R} \left| B_R^*(A_Q^*(f)) \right|^{q_1}
d\mu \right)^{1/q_1 } \leq c \|f\|_{BMO_{q_0}}.
%\tag{$A_0$}
\end{equation}
In some cases, we will require the following stronger assumption,
which describes that the operators $A_Q^*$ continuously act on the
$Bmo$ spaces~: there is a constant $c$ such that
\begin{equation} \label{ass:cont2}
\|A_Q^*(f)\|_{BMO_{q_1}} \leq c \|f\|_{BMO_{q_0}}.
%\tag{$A$}
\end{equation}

\gb Then we have the two following results.

\begin{thm} \label{thm} Let us first assume that the operators $B_Q$ depend only on the radius $r_Q$ of the ball and that (\ref{assum}), (\ref{assum2}) and (\ref{ass:cont}) hold. Then the spaces $Bmo_{q_0}$ satisfies to John-Nirenberg inequalities. There exist constants $\rho_1,\rho_2>0$ such that for all function $f\in Bmo_{q_0}$ and every ball $Q\subset X$
$$ \mu\left(\left\{x\in Q,\ \left|B_Q^*(f)\right| >\lambda \|f\|_{BMO_{q_0}} \right\}\right) \leq \rho_1\mu(Q) \left[e^{-\rho_2\lambda} +\frac{{\bf 1}_{q_1<\infty}}{\lambda^{q_1}}\right].$$
\end{thm}

\begin{cor} \label{cor:imp} If we are working on the Euclidean space $\R^n$, we can just require
\begin{equation} \label{ass:cont3}
\sup_{\genfrac{}{}{0pt}{}{R \textrm{ ball}}{R\subset Q}} \left(\frac{1}{\mu(R)}\int_{R} \left| B_R^*(A_Q^*(f)) \right|^{q_1} d\mu \right)^{1/q_1 } \leq c \|f\|_{BMO_{q_0}}. %\tag{$A_0'$}
\end{equation}
instead of (\ref{ass:cont}), see Remark \ref{rem:imp}.
\end{cor}

\mb
Let us now consider general operators $B_Q$, which could depend on the ball $Q$.

\begin{thm} \label{thm2} Let $\epsilon >0$ and suppose that the coefficients $(\gamma_j)_j$ given by (\ref{assum}) satisfy
$$\sup_{j} \gamma_j <\infty$$
instead of the stronger inequality (\ref{assum2}). \\
Under (\ref{ass:cont2}), the spaces $Bmo_{\epsilon,q_0}$ satisfies to John-Nirenberg inequalities. There exist constants $\rho_1,\rho_2>0$ such that for all function $f\in Bmo_{\epsilon,q_0}$ every ball $Q\subset X$ and index $k\geq 0$
\be{eq:dfg} \mu\left(\left\{x\in C_k(Q),\ \left|B_Q(f)\right| >\lambda \|f\|_{BMO_{\epsilon,q_0}} \right\}\right) \leq \rho_1 2^{-\epsilon k} \mu(2^k Q) \left[e^{-\rho_2\lambda} +\frac{{\bf 1}_{q_1<\infty}}{\lambda^{q_1}}\right]. \ee
\end{thm}

\begin{rem} \label{rem:gamma} Indeed, we prove a more accurate result in the particular case of non-increasing coefficients $\gamma_j$. In this case, we can work with the smaller norm
 $$ \left\|f \right\|_{\widetilde{BMO}_{\epsilon,q}} := \sup_{Q \textrm{ ball}}\ \sup_{\genfrac{}{}{0pt}{}{j\geq 0}{\gamma_j\neq 0}}\ 2^{j\epsilon}\left(\frac{1}{\mu(2^j Q)}\int_{C_j(Q)} \left| B_Q^*(f) \right|^q d\mu \right)^{1/q }$$
and prove (\ref{eq:dfg}). We let the details to the reader as it just suffices to follow the contribution of coefficients $\gamma_j$ in the proof.
\end{rem}

\mb As usual, a John-Nirenberg inequality permits to prove some equivalence in $BMO$ spaces, making varying the exponent.

\begin{cor} \label{corco}
Under the assumptions of Theorem \ref{thm} (resp. Theorem \ref{thm2}), the norms $BMO_{q}$ (resp. $BMO_{\epsilon,q}$ for some $\epsilon>0$) for $q\in(q_0,q_1)$ are equivalent and consequently the spaces $BMO_q$ are equal.
\end{cor}

\dem  We only treat the case of the $BMO_{q}$ spaces as it is exactly the same proof for the spaces $BMO_{\epsilon,q}$, using Theorem \ref{thm2} instead of Theorem \ref{thm}. We take two exponents $ r> q$ belonging to the range $(q_0,q_1)$ and a function
$ \phi \in Bmo_r\cap Bmo_q$. First, using H\"older inequality, we have for each ball $Q$
$$ \left(\frac{1}{\mu(Q)}\int_{Q} \left| B_Q^*(\phi) \right|^q d\mu
\right)^{1/q } \leq  \left( \frac{1}{\mu(Q)}
\int_{Q} \left| B_Q^*(\phi) \right|^r d\mu \right)^{1/r }
\leq  \| \phi \|_{BMO_r},$$ therefore we have
the inclusion $ Bmo_r \subset Bmo_q.$ \\
Then it remains to check that  $ Bmo_q \subset Bmo_r.$ Using John-Nirenberg inequality (obtained in Theorem \ref{thm}), we get a weak inequality for every ball $Q$
$$ \frac{1}{\mu(Q)^{1/q_1}}\left\| B_Q^*(\phi) \right\|_{L^{q_1,\infty}} \lesssim \| \phi \|_{BMO_{q_0}} \lesssim \| \phi \|_{BMO_{q}}.$$
We conclude by invoking Kolmogorov's inequality to get
$$ \left( \frac{1}{\mu(Q)}
\int_{Q} \left| B_Q^*(\phi) \right|^r d\mu \right)^{1/r } \leq \left( \frac{q_1}{q_1-r}\right)^{1/r} \frac{1}{\mu(Q)^{1/q_1}}\left\| B_Q^*(\phi) \right\|_{L^{q_1,\infty}},$$
which finally yields
$$\left( \frac{1}{\mu(Q)}
\int_{Q} \left| B_Q^*(\phi) \right|^r d\mu \right)^{1/r } \lesssim \| \phi \|_{BMO_{q}}.$$
\findem

\subsection{Proof of Theorems \ref{thm} and \ref{thm2}} \label{subsec:proof}

\mb The following proof has been written with abstract operators $B_Q$ (depending on the ball), as we will refer to it for Theorem \ref{thm2} requiring this abstract framework.

 {\noindent {\bf Proof of Theorem \ref{thm}: }} We follow the ideas of \cite{DY}. \\
By homogeneity, we can assume that $\|f\|_{BMO_{q_0}}=1$, so we have to prove for any fixed ball $Q$
\be{eq:amontrer} \mu\left(\left\{x\in Q,\ \left|B_Q^*(f)\right| >\lambda \right\}\right) \leq \rho_1\mu(Q) \left[e^{-\rho_2\lambda} +\frac{{\bf 1}_{q_1<\infty}}{\lambda^{q_1}}\right]. \ee
Obviously (\ref{eq:amontrer}) holds for $\lambda\leq 1$ with $\rho_2=1$ and $\rho_1=e$. So we only consider $\lambda>1$ and set
$$ f_0:={\bf 1}_{Q} B_Q^*(f).$$
Then we get
$$ \|f_0\|_{L^1} \leq \int_{Q} |B_Q^*(f)| d\mu \leq \|f\|_{BMO_{q_0}} \mu(Q) \leq \mu(Q).$$
Using a constant $\beta>1$ (later choosen) and the Hardy-Littlewood maximal operator $\M$, we set
$$ F:=\{x,\ \M(f_0)(x) \leq \beta \} \qquad \Omega:=F^c = \{x,\ \M(f_0)(x) >\beta\}.$$
We take a Whitney decomposition of $\Omega$ : that is a collection of balls $(Q_{1,i})_{i}$ such that
\begin{itemize}
 \item[ ] \hspace{1cm} $a-) \quad$  $\Omega=\bigcup_{i} Q_{1,i}$
 \item[ ] \hspace{1cm} $b-) \quad$ each point is contained in at most a finite number $N$ of balls $Q_{1,i}$
$$ \sum_{i} {\bf 1}_{Q_{1,i}} \leq N$$
 \item[ ] \hspace{1cm} $c-) \quad$ there exists $\kappa>1$ such that for all $i$, $\kappa Q_{1,i} \cap F \neq \emptyset$.
\end{itemize}
From $a-)$, for all $x\in Q\setminus (\cup_i Q_{1,i})$
\be{eq:pre} |B_Q^*(f)(x)| = |f_0(x)| \leq \M(f_0)(x) \leq \beta. \ee
The weak-type $(1,1)$ of the Hardy-Littlewood maximal operator yields
\be{eq:1} \sum_{i} \mu(Q_{1,i}) \leq N\mu(\Omega) \lesssim \frac{1}{\beta}\|f_0\|_{L^1} \leq \frac{c_1}{\beta}\mu(Q) \ee
for some numerical constant $c_1>0$. \\
We choose $\beta$ such that for all $Q_{1,i} \cap Q \neq \emptyset$,
$Q_{1,i}\subset 2 Q$. Let us check that this is possible. Indeed for
such balls, if $r_{Q_{1,i}}\leq r_Q$ then we have nothing to do,
else we have from (\ref{eq:1})~:
$$ \mu(Q) \lesssim \left(\frac{r_Q}{r_{Q_{1,i}}}\right)^{\delta} \mu(Q_{1,i}) \leq \left(\frac{r_Q}{r_{Q_{1,i}}}\right)^{\delta} \frac{c_3}{\beta} \mu(Q).$$
So we choose $\beta$ is order that the previous inequality implies $r_{Q} \geq 2 r_{Q_{1,i}}$, which yields $Q_{1,i}\subset 2 Q$.

\mb Then we are interested in the following Lemma.
\begin{lem} \label{lem:lem} There exists a numerical constant $c_2\geq 1$, such that for all $i$~:
\be{eq:est2} \left(\frac{1}{\mu(Q_{1,i})}\int_{Q_{1,i}} \left|B_{Q}^*(f)-B_{Q_{1,i}}^*(f) \right|^{q_1} d\mu \right)^{1/q_1} \leq c_2 \beta. \ee
\end{lem}

\mb For readibility, we postpone the proof. Let now come back to the proof of our main Theorem.
For each $i$, we repeat the procedure as follows~: consider
$$ f_{1,i} := {\bf 1}_{Q_{1,i}} B_{Q_{1,i}}^*(f).$$
Then we consider a collection of balls $(Q_{2,i,m})_{m}$ such that
\begin{itemize}
 \item for all $x\in Q_{1,i}\setminus \left(\cup_{m} Q_{2,i,m}\right)$
$$ \left|B_{Q_{1,i}}^*(f)(x)\right|\leq \beta$$
 \item we have
$$ \sum_{m} \mu(Q_{2,i,m}) \leq \frac{c_1}{\beta}\mu(Q_{1,i})$$
\item for all balls $Q_{2,i,m}$ intersecting $Q_{1,i}$, we have
$$ \left(\frac{1}{\mu(Q_{2,i,m})}\int_{Q_{2,i,m}} \left|B_{Q_{1,i}}^*(f)-B_{Q_{2,i,m}}^*(f) \right|^{q_1} d\mu \right)^{1/q_1} \leq c_2 \beta.$$
\end{itemize}
Then we put together all families $(Q_{2,i,m})_{m}$ for all indices $i$ and we get a new family $(Q_{2,m})_{m}$.
We check that
$$ \sum_{m} \mu(Q_{2,m}) \leq \frac{c_1}{\beta} \sum_{i} \mu(Q_{1,i}) \leq \left(\frac{c_1}{\beta}\right)^2 \mu(Q).$$
Moreover for all $x\in Q\setminus \left(\cup_{i} Q_{1,i}\right)$, we already know from (\ref{eq:pre})
$$ \left|B_Q^*(f)(x)\right| \leq \beta.$$
For all $x$ belonging to one ball $Q\cap Q_{1,i}$ but not in the associated collection $(Q_{2,i,m})_{m}$, we have
$$ \left|B_{Q_{1,i}}^*(f)(x)\right| \leq \beta,$$
so
$$\left|B_Q^*(f)(x)\right| \leq \beta + \left| B_{Q}^*(f)(x)-B_{Q_{1,i}}^*(f)(x)\right|.$$
According to (\ref{eq:est2}) we get~:
$$ \left(\frac{1}{\mu(Q_{1,i})}\int_{Q_{1,i}\setminus \left(\cup_m Q_{2,i,m}\right) } \left|B_{Q}^*(f)\right|^{q_1} d\mu \right)^{1/q_1} \leq (c_2+1) \beta \leq 2c_2\beta.$$

\gb We iterate this procedure and then associated to a collection $(Q_{k,i})_{i}$, we build for all $i$ a collection $(Q_{k+1,i,m})_{m}$ and also a collection $(Q_{k+1,m})_{m}=\cup_{i} (Q_{k+1,i,m})_{m}$ satisfying~:
\begin{itemize}
 \item for all $x\in Q_{k,i}\setminus \left(\cup_{m} Q_{k+1,i,m}\right)$
\be{eq:azebis} \left|B_{Q_{k,i}}^*(f)(x)\right|\leq \beta \ee
 \item we have
$$ \sum_{m} \mu(Q_{k+1,i,m}) \leq \frac{c_1}{\beta}\mu(Q_{k,i})$$
\item for all balls $Q_{k+1,i,m}$ intersecting $B_{k,i}$, we have
\be{eq:aze2} \left(\frac{1}{\mu(Q_{k+1,i,m})}\int_{Q_{k+1,i,m}} \left|B_{Q_{k,i}}^*(f)-B_{Q_{k+1,i,m}}^*(f) \right|^{q_1} d\mu \right)^{1/q_1} \leq c_2 \beta.\ee
\end{itemize}
So it results that for all integer $k$
\be{eq:res1} \sum_{m} \mu(Q_{k,m}) \leq \frac{c_1}{\beta} \sum_{i} \mu(Q_{k-1,i}) \leq \left(\frac{c_1}{\beta}\right)^{k} \mu(Q).\ee

\gb {\bf First case:} If $q_1<\infty$. \\
For $\lambda$ large enough (we have seen at the beginning of the proof that (\ref{eq:amontrer}) holds for $\lambda \lesssim 1$), we denote a positive integer $K\geq K_0$ and a constant $\gamma<1$ such that
$$ \gamma^K\lambda\simeq 2\beta \textrm{  and  } \left(\frac{c_1}{\beta}\right)^K\leq \lambda^{-q_1}. $$
It is possible with $\gamma^{q_1}=lc_1/\beta$ and $l\geq (2\beta)^{q_1/K_0}$ ($\beta>c_1$ beeing chosen large enough, we can found an integer $K_0$ such that $l=(2\beta)^{q_1/K_0}>1$ satisfies $lc_1/\beta>1$). \\
We remark that we have in particular
\be{eq:aze4} \frac{1}{\gamma^{q_1}} \frac{c_1}{\beta} \leq 1. \ee
Now we can estimate as follows
\begin{align*}
\mu\left(\left\{x\in Q,\ |B_{Q}^*(f)(x)|>\lambda\right\}\right) & \leq \sum_{i} \mu\left(\left\{x\in Q_{1,i},\ |B_{Q}^*(f)(x)|>\lambda\right\}\right) \\
& \leq \sum_{i} \mu\left(\left\{x\in Q_{1,i},\ |B_{Q}^*(f)(x)-B_{Q_{1,i}}^*(f)(x)|>(1-\gamma)\lambda\right\}\right) \\
& \hspace{0.5cm} + \sum_{i} \mu\left(\left\{x\in Q_{1,i},\ |B_{Q_{1,i}}^*(f)(x)|>\gamma\lambda\right\}\right).
\end{align*}
The first term is bounded by
\begin{align*}
\sum_i \mu\left(\left\{x\in Q_{1,i},\ |B_{Q}^*(f)(x)-B_{Q_{1,i}}^*(f)(x)|>(1-\gamma)\lambda\right\}\right) & \\
& \hspace{-3cm} \leq (1-\gamma)^{-q_1}\lambda^{-q_1} \sum_i \int_{Q_{1,i}} |B_{Q}^*(f)-B_{Q_{1,i}}^*(f)|^{q_1} d\mu  \\
& \hspace{-3cm} \leq \left(\frac{c_2\beta}{(1-\gamma)\lambda}\right)^{q_1} \sum_i \mu(Q_{1,i})  \\
& \hspace{-3cm} \leq \left(\frac{c_2\beta}{(1-\gamma)\lambda}\right)^{q_1} \frac{c_1}{\beta} \mu(Q).
\end{align*}
Then we repeat the procedure with $B_{Q_{1,i}}^*(f)$ instead of $B_{Q}^*(f)$
\begin{align*}
\sum_{i} \mu\left(\left\{x\in Q_{1,i},\ |B_{Q_{1,i}}^*(f)(x)|>\gamma\lambda\right\}\right) & \\
& \hspace{-3cm} \leq \sum_{i,j} \mu\left(\left\{x\in Q_{1,i,j},\ |B_{Q_{1,i}}^*(f)(x)|>\gamma\lambda\right\}\right) \\
& \hspace{-3cm} \leq \sum_{i,j} \mu\left(\left\{x\in Q_{1,i,j},\ |B_{Q_{1,i}}^*(f)(x)-B_{Q_{1,i,j}}^*(f)(x)|>(1-\gamma)\gamma\lambda\right\}\right) \\
& \hspace{-2.5cm} + \sum_{i,j} \mu\left(\left\{x\in Q_{1,i,j},\ |B_{Q_{1,i,j}}^*(f)(x)|>\gamma^2\lambda\right\}\right).
\end{align*}
The first term (in the last inequality) is now controled by
\begin{align*}
\sum_{i,j} \mu\left(\left\{x\in Q_{1,i,j},\ |B_{Q_{1,i}}^*(f)(x)-B_{Q_{1,i,j}}^*(f)(x)|>(1-\gamma)\gamma\lambda\right\}\right) \\
& \hspace{-5cm} \leq [(1-\gamma)\gamma \lambda]^{-q_1} \sum_{i,j} \int_{Q_{1,i,j}} |B_{Q_{1,i}}^*(f)-B_{Q_{1,i,j}}^*(f)|^{q_1} d\mu \\
& \hspace{-5cm} \leq \left(\frac{c_2\beta}{(1-\gamma)\gamma\lambda}\right)^{q_1} \sum_{i,j} \mu(Q_{1,i,j}) \\
& \hspace{-5cm} \leq \left(\frac{c_2\beta}{(1-\gamma)\gamma\lambda}\right)^{q_1} \left(\frac{c_1}{\beta}\right)^2 \mu(Q)
\end{align*}
and the second one is equal to
$$ \sum_{i} \mu\left(\left\{x\in Q_{2,i},\ |B_{Q_{2,i}}^*(f)(x)|>\gamma^2\lambda \right\}\right).$$
Thus, by iterating this reasoning, we deduce that
\begin{align*}
\mu\left(\left\{x\in Q,\ |B_{Q}^*(f)(x)|>\lambda\right\}\right) & \leq \frac{\mu(Q)}{\lambda^{q_1}} \sum_{k=0}^{K-1} \left(\frac{c_2\beta}{(1-\gamma)\gamma^k}\right)^{q_1} \left(\frac{c_1}{\beta}\right)^{k+1} \\
 & \ + \sum_{i} \mu\left(\left\{x\in Q_{K,i},\ |B_{Q_{K,i}}^*(f)(x)|>\gamma^K\lambda\right\}\right).
\end{align*}
Consequently, since (\ref{eq:res1})
\begin{align*}
\mu\left(\left\{x\in Q,\ |B_{Q}^*(f)(x)|>\lambda\right\}\right) & \leq \frac{\mu(Q)}{\lambda^{q_1}} \sum_{k=0}^{K-1} \left(\frac{c_2\beta}{(1-\gamma)\gamma^k}\right)^{q_1} \left(\frac{c_1}{\beta}\right)^{k+1} \\
 & \ + \left(\frac{c_1}{\beta}\right)^{K}\mu(Q).
\end{align*}
By the choice of the constant $\gamma$ and the integer $K$ and (\ref{eq:aze4}), it comes
$$ \mu\left(\left\{x\in Q,\ |B_{Q}^*(f)(x)|>\lambda\right\}\right) \lesssim \lambda^{-q_1} \mu(Q)$$
which corresponds to the desired inequality (\ref{eq:amontrer}) when $q_1<\infty$.

\gb {\bf Second case:} If $q_1=\infty$.\\
In this case, we repeat the proof of \cite{DY}. For $\lambda$ large enough, we denote $K\geq K_0$ a lare enough integer such that
$$ Kc_2 \beta < \lambda \leq (K+1)c_2 \beta. $$
Then, since (\ref{eq:est2}) and (\ref{eq:aze2}), it follows that on $Q_{1,i}\setminus (\cup_j Q_{1,i,j})$, $|B_Q^*(f)|\leq \beta$, on $Q_{1,i_1} \cap Q_{2,i_2}\setminus (\cup_j Q_{2,i_2,j})$
$$ |B_Q^*(f)|\leq |B_Q^*(f)-B_{Q_{1,i}}^*(f)| + |B_{Q_{1,i_1}}^*(f)| \leq (1+c_2)\beta \leq 2c_2\beta$$
and by iterating on $Q_{1,i_1} \cap \cdots \cap Q_{K-1,i_{K-1}}\setminus (\cup_j Q_{K,i_{K-1},j})$
\begin{align*}
 |B_Q^*(f)|& \leq |B_Q^*(f)-B_{Q_{1,i}}^*(f)| + \sum_{l=1}^{K-2} |B_{Q_{l+1,i_l}}^*(f)-B_{Q_{l+1,i_{l+1}}}^*(f)| + |B_{Q_{K-1,i_{K-1}}}^*(f)| \\
& \leq (1+(K-1) c_2)\beta \leq Kc_2\beta<\lambda.
\end{align*}
Hence
$$ \left\{x\in Q,\ |B_{Q}^*(f)(x)|>\lambda\right\} \subset \bigcup_{i} Q_{K,i},$$
which yields thanks to (\ref{eq:res1})
$$ \mu\left(\left\{x\in Q,\ |B_{Q}^*(f)(x)|>\lambda\right\}\right) \leq \left(\frac{c_1}{\beta}\right)^{K}.$$
As $c_1/\beta$ is a constant smaller than $1$ and $K\simeq \lambda$, this allows us to obtain the desired inequality.
\findem

\mb It remains us to prove Lemma \ref{lem:lem}. We recall the statement with the notations of the previous proof.
\begin{lem}  \label{lem:lem2} There exists a numerical constant $c_2\geq 1$, such that for all $i$~:
\be{eq:est2bis} \left(\frac{1}{\mu(Q_{1,i})}\int_{Q_{1,i}} \left|B_{Q}^*(f)-B_{Q_{1,i}}^*(f) \right|^{q_1} d\mu \right)^{1/q_1} \leq c_2 \beta. \ee
\end{lem}

\dem The desired result corresponds to a ``local version'' of inequality (3.3) in \cite{DY} (extended in our abstract framework), which essentialy rests on Proposition 2.6 of \cite{DY}. We know that $Q_{1,i}\subset 2Q$ and we have
\begin{align}
B_{Q}^*(f)-B_{Q_{1,i}}^*(f) & = A_{Q}^*(f)-A_{Q_{1,i}}^*(f) = \left[A_Q^*(f) - A_{Q_{1,i}}^*A_Q^*(f) \right]+\left[A_{Q_{1,i}}^*A_Q^*(f) - A_{Q_{1,i}}^*(f)\right] \nonumber \\
 & = B_{Q_{1,i}}^*A_Q^*(f) - A_{Q_{1,i}}^*B_Q^*(f). \label{eq:estt}
\end{align}
Let us study the first term. As $Q_{1,i}\subset 2Q$, Assumption (\ref{ass:cont}) implies~:
$$\left(\frac{1}{\mu(Q_{1,i})}\int_{Q_{1,i}} \left|B_{Q_{1,i}}^*A_Q^*(f)(f) \right|^{q_1} d\mu \right)^{1/q_1} \lesssim \|f\|_{BMO_{q_0}} \lesssim 1\leq \beta,$$
as $\beta$ is chosen large enough.\\
So it remains to study the second term $A_{Q_{1,i}}^*B_Q^*(f)$. To estimate it, we have to use $(L^{q_0}-L^{q_1})$-off diagonal decays of $A_{Q_{1,i}}$ as follows. Let denote $P$ the first integer such that $2Q\subset 2^{P+1}Q_{1,i}$ and $2Q \cap (2^{P}Q_{1,i})^c\neq \emptyset$. Then since Assumption (\ref{assum}), we decompose
$$
 \left(\frac{1}{\mu(Q_{1,i})}\int_{Q_{1,i}} \left|A_{Q_{1,i}}^* B_{Q}^*(f) \right|^{q_1} d\mu \right)^{1/q_1} \lesssim I+II
$$
with
$$I:= \sum_{j=0}^{P+1} \gamma_{j} \left(\frac{1}{\mu(2^jQ_{1,i})}\int_{C_j(Q_{1,i})} \left| B_{Q}^*(f) \right|^{q_0} d\mu \right)^{1/q_0}$$
and
$$II:=\sum_{j=P+2}^{\infty} \gamma_{j} \left(\frac{1}{\mu(2^jQ_{1,i})}\int_{C_j(Q_{1,i})} \left| B_{Q}^*(f) \right|^{q_0} d\mu \right)^{1/q_0}.$$
It follows from the property of ball $Q_{1,i}$ (property $c-)$ in the proof of Theorem \ref{thm}), that there exists another constant $\kappa'$ with for $j\leq P+1$
$$ \left(\frac{1}{\mu(2^jQ_{1,i})}\int_{C_j(Q_{1,i})} \left| B_{Q}^*(f) \right|^{q_0} d\mu \right)^{1/q_0} \leq \kappa' \left(\frac{1}{\mu(2^j \kappa Q_{1,i})}\int_{2^j \kappa Q_{1,i}} \left| f_0 \right|^{q_0} d\mu \right)^{1/q_0} \leq \kappa' \beta.$$
So it yields
$$ I\leq  \sum_{j=0}^{P+1} \gamma_{j} \kappa' \beta \lesssim \beta.$$
To estimate the second term $II$. For any $j\geq P+1$, we know that $2^{j}Q_{1,i}$ contains the ball $Q$ and so $2^jr_{Q_{1,i}} \geq r_{Q}$. So we choose $(\tilde{Q}^j_k)_{k}$ a bounded covering of $C_j(Q_{1,i})$ with balls of radius $r_{Q}$ and as previously, we get (using that the operators $B_Q$ only depend on the radius of the ball $Q$)
\begin{align}
II & \leq \sum_{j=P+2}^{\infty} \gamma_{j} \left(\frac{1}{\mu(2^j Q_{1,i})}\sum_{k} \int_{\tilde{Q}^j_k} \left| B_{Q}^*(f) \right|^{q_0} d\mu \right)^{1/q_0} \nonumber \\
& \lesssim \sum_{j=P+2}^{\infty} \gamma_{j} \left(\frac{1}{\mu(2^j Q_{1,i})}\sum_{k} \int_{\tilde{Q}^j_k} \left| B_{\tilde{Q}^j_k}^*(f) \right|^{q_0} d\mu \right)^{1/q_0} \label{eq:aze} \\
& \lesssim \beta \sum_{j=P+2}^{\infty} \gamma_{j} \left(\frac{1}{\mu(2^j Q_{1,i})}\sum_{k} \mu(\tilde{Q}^j_k) \right)^{1/q_0} \nonumber \\
& \lesssim \beta \sum_{j=P+2}^{\infty} \gamma_{j} \lesssim \beta. \nonumber
\end{align}
So finally, the estimates of $I$ and $II$ imply (\ref{eq:est2bis}), which concludes the proof.
\findem

\begin{rem} \label{rem:imp}
Let us show how can we obtain Corollary \ref{cor:imp}. As explained in the proof of Theorem \ref{thm}, we use a Whithney decomposition of the set
$$ \Omega:= \{x,\ \M\left({\bf 1}_{Q} B_Q^*(f) \right)(x) >\beta\},$$
for $Q$ a ball and $\beta>1$ some fixed parameter.
Using the dyadique structure of $\R^d$, let us deal we a dyadic cube $Q$. We can choose a Whithney decomposition of $\Omega$ with dyadic (relatively to $Q$) cubes $Q_{1,i}$ -- see Theorem 5.2 of \cite{JM} for a detailed construction --. Then the proof is based on such balls $Q_{1,i}$ such that
$$ Q_{1,i} \cap Q \neq \emptyset , \quad \textrm{and} \quad Q_{1,i} \subset cQ,$$
for some constant $c>1$. We have chosen $c=2$ for simplicity but we
can consider $c=3/2$ for example. Then the dyadic structure of the
Euclidean space, implies that $Q_{1,i}$ is included in $Q$. Then we
reproduce the same argument and Assumption (\ref{ass:cont3}) is
sufficient to conclude.
\end{rem}

\gb
Let us now consider general operators $B_Q$, which could depend on the ball $Q$.

\mb To get a result concerning abstract operators $B_Q$ (they depend now on the ball and not only on the radius), we have to require some extra properties. In the previous proof, the only one point where we used the property of dependence (on the radii) of the operators $B_Q$, is the inequality (\ref{eq:aze}). So let us just take the notations of the previous proof and recall the problem: for $j\geq P+1$, $C_j(Q_{1,i}) \subset C_{j-P}(Q)$ and we have to estimate
$$\int_{C_j(Q_{1,i})} \left| B_{Q}^*(f) \right|^{q_0} d\mu.$$
From the $BMO$-norm, we only have information about $B_Q^*(f)$ on the ball $Q$ so we do not know how can we control this term. In order to get around this lack of information, we use the $BMO_{\epsilon}$ associated to the notion of {\it molecules} (see Theorem \ref{thm2}).

{\noindent {\bf Proof of Theorem \ref{thm2}: }} For $k=0$, we follow the proof of Theorem \ref{thm}. In this one, the only one point where we used the property of dependence (on the radii) of the operators $B_Q$, is the inequality (\ref{eq:aze}). So let us just take the notations of the previous proof and check how can we get around this problem. So we have an integer $j\geq P+1$ (so $C_j(Q_{1,i}) \subset C_{j-P}(Q)$) and we have to estimate
$$\int_{C_j(Q_{1,i})} \left| B_{Q}^*(f) \right|^{q_0} d\mu.$$
We also have
$$ \int_{C_j(Q_{1,i})} \left| B_{Q}^*(f) \right|^{q_0} d\mu \leq \int_{C_{j-P}(Q)} \left| B_{Q}^*(f) \right|^{q_0} d\mu \lesssim 2^{-\epsilon(j-P)}.$$
Consequently, we get
\begin{align*}
II & \leq \sum_{j=P+1}^{\infty} \gamma_{j} 2^{-\epsilon(j-P)} \nonumber \\
& \lesssim \sum_{j=P+2}^{\infty} 2^{-\epsilon(j-P)} \lesssim 1 \lesssim \beta.
\end{align*}
This estimate permits to conclude the proof of Lemma \ref{lem:lem2} and by this way the proof of the desired inequality for $k=0$.\\
For $k\geq 1$, we produce the same reasoning, starting from the function
$f_0:={\bf 1}_{C_k(Q)} B_Q^*(f)$ which satisfies
$$ \|f_0\|_{L^1} \leq \int_{C_k(Q)} |B_Q^*(f)| d\mu \leq \|f\|_{BMO_{\epsilon,q_0}} 2^{-\epsilon k} \mu(2^k Q).$$
We reproduce the same proof (using Assumption (\ref{ass:cont2})), which we left to the reader.
\findem

\section{Check-out of our assumptions in some usual cases} \label{sec:check}

We devote this section to check that our results generalize so already known particular cases and more precisely that our assumptions are satisfying.

\subsection{The John-Nirenberg space}

Consider the Euclidean space $ X=\R^n$ and the usual BMO space. In \cite{JN}, the first John-Nirenberg inequalities was proved using a Calder\'on-Zygmund decomposition (similar to what we previously done). We refer the reader to Chap IV 1.3 of the book \cite{stein} for another proof based on the duality $H^1-BMO$.

\mb This first BMO space is defined by the operators~:
$$ B_Q^*(f):= f - \left(\frac{1}{\mu(Q)}\int_Q f \right){\bf 1}_{Q}.$$
So $A_Q^*$ is the mean value operator and obviously off diagonal decays (\ref{assum}) hold with $q_1=\infty$ and $q_0=1$.
In this case Assumption (\ref{ass:cont2}) does not hold. However for such operators, the coefficients $\gamma_{j}=0$ as soon as $j\geq 1$.
So indeed to apply Theorem \ref{thm2}, we have just to check Assumption (\ref{ass:cont3}), thanks to Corollary \ref{cor:imp} and Remark \ref{rem:gamma} (in this particular case the spaces $\widetilde{BMO}_{\epsilon,q}$ are equal to $BMO_{q}$ because $\gamma_j=0$ for $j\geq 1$), which is
%So indeed to apply Theorem \ref{thm2}, we have just to check Assumption (\ref{ass:cont}), unfortunately does not hold too.
%Since Remark \ref{rem:imp}, as we are working on the Euclidean space, it is sufficient to check Assumption (\ref{ass:cont3}), which is
$$ \sup_{\genfrac{}{}{0pt}{}{R \textrm{ ball}}{R\subset Q}} \left(\frac{1}{\mu(R)}\int_{R} \left| B_R^*(A_Q^*(f)) \right|^{q_1} d\mu \right)^{1/q_1 } \leq c \|f\|_{BMO_{q_0}}. $$
And this one is satisfied as for $R\subset Q$, we have~:
$$ {\bf 1}_{R} B_R^*(A_Q^*(f)) = \left(\frac{1}{\mu(Q)}\int_Q f \right) \left[ {\bf 1}_{R} - {\bf 1}_{R}\right] = 0.$$

\gb
{\bf Conclusion :} In the framework of the classical BMO space, our Assumptions (\ref{assum}) and (\ref{ass:cont3}) are satisfied. We recover also the John-Nirenberg inequality (see \cite{JN}).

\subsection{The Morrey-Campanato spaces}

Consider the set $ X=[0,1]$ with its Euclidean structure. We refer the reader to \cite{Li} and \cite{DDY} for works related to Morrey-Campanato spaces and associated John-Nirenberg inequalities.

\mb Let us first define these spaces.

\begin{df} For $\beta\geq 0$, $s\in {\mathbb N}$ and $q\in(1,\infty)$, we say that a locally integrable function $f\in L^1(X)$ belongs to the Morrey-Campanato spaces $L(\beta,q,s)$ if
$$ \|f\|_{L(\beta,q,s)}:=\sup_{Q\in {\mathcal Q}} |Q|^{-\beta} \left[\aver{Q} \left|f(x)-P_Q(f)(x) \right|^q dx \right]^{1/q} <\infty,$$
where for $Q$ a ball (an interval) of $X$, $P_Q(f)$ is the unique polynomial function of degree at most $s$ such that for all $i\in\{0,..,s\}$
$$ \int_Q x^i \left( f(x)-P_Q(f)(x)\right) dx = 0.$$
\end{df}

\begin{rem} So we remark that $L(\beta,q,0)$ exactly corresponds to the previous BMO space as in this case $P_Q(f) = \aver{Q}f$.
\end{rem}

\mb In this framework, we set $A_Q:=P_Q^{*}$ in order that $L(0,q,s)$ can be identified to our BMO space. An easy computation gives that $P_Q(f)$ is a polynomial function whose coefficients are given by the quantities
$$ \int_Q f(x) x^{i} dx$$
for $i\in\{0,..,s\}$. So $A_Q^*=P_Q$ can be written as follows
$$A_Q^*(f)(x) = \sum_{j=0}^s c_j x^{j} {\bf 1}_{Q}(x)$$
with coefficients $c_j$ satisfying
$$ \left|c_j\right| \lesssim \int_Q|f(x)|dx$$
since we are working on $X=[0,1]$. It comes that off- diagonal
decays (\ref{assum}) hold with $q_1=\infty$ and $q_0=1$. As
previously since we are working on the Euclidean space and
coefficients $\gamma_j=0$ as soon as $j\geq 1$, it is sufficient to
check~:
$$ \sup_{\genfrac{}{}{0pt}{}{R \textrm{ ball}}{R\subset Q}} \left(\frac{1}{\mu(R)}\int_{R} \left| B_R^*(A_Q^*(f)) \right|^{q_1} d\mu \right)^{1/q_1 } \leq c \|f\|_{BMO_{q_0}}. $$
And this property is satisfied as for $R\subset Q$, we have~:
$$ B_R^*(A_Q^*(f)) = P_Q(f)-P_RP_Q(f) = 0.$$
The last equality is due to the fact that $P_Q(f)$ is a polynomial function of degree at most $s$ so by uniqueness (in the definition of $P_R$): $P_R[P_Q(f)]=P_Q(f)$.

\gb
{\bf Conclusion :} In the framework of the classical Morrey-Campanato spaces $L(0,q,s)$, our Assumptions (\ref{assum}) and (\ref{ass:cont3}) are satisfied. We recover also the John-Nirenberg inequality for all $q\in(1,\infty)$ (see \cite{Li}). For $\beta>0$, we refer the reader to a forthcoming work of the first author and J.M. Martell \cite{BM} (and \cite{JM,J}), dealing with more general self-improvement properties of inequalities.

% \subsection{BMO spaces for Schr\"odinger operators with potentials.}
%
% Consider the Euclidean space $X=\R^n$ and the framework of \cite{DGMTZ}. J.~Dziubanski1, G.~Garrig\'os, T.~Mart\'inez, J.L.~Torrea, J.~Zienkiewicz consider a non-negative potential $V$ belonging to some Reverse Holder class $RH_s$ with $s>n/2$. \\
% By this assumption, it is well known that $-L=-\Delta+V$ generates a
% $L^2$-bounded semigroup $(K_t)_{t>0}$, whose kernels satisfy some
% gaussian estimates. J.~Dziuba\'nski and J.~Zienkiewicz studied an associated Hardy space $H^1_L$.
% Then in \cite{DGMTZ}, the dual space $BMO$ is considered. It corresponds to the following particular choice of operators $B_Q$.
% First, we set
% $$ \rho(x):= \sup \left\{ r>0,\ \frac{1}{r^{n-2} \int_{B(x,r) V(y) dy \leq 1 \right\}.$$
% For all ball $Q$, we define
% $$ B_Q^*(f)(x) := \left\{ \begin{array}{l}
%  f(x){\bf 1}_{Q}(x)  \qquad \textrm{  if $r_Q\geq \rho(x)$} \\
%  f(x){\bf 1}_{Q}(x)-\left(\frac{1}{|Q|} \int_Q f(y)dy \right){\bf 1}_{Q}(x) \qquad  \textrm{  if $r_Q\leq \rho(x)$}
% \end{array}
% \right. . $$
% Then our BMO spaces correspond to the ones defined in \cite{DGMTZ}.
%
% \mb As previously, Assumption (\ref{assum}) holds for $q_1=\infty$ and $q_0=1$.

\subsection{General case of semigroup}

Let us recall the framework of \cite{DY}.

\mb Consider a space of homogeneous type $(X,d,\mu)$ with a family of operators $({\mathcal A}_r)_{r>0}$ satisfying according to a parameter $m>0$ the following properties~:
\begin{itemize}
 \item For every $r>0$, the linear operator ${\mathcal A}_r$ is given by a kernel $a_r$ satisfying
$$ \left|a_{r}(x,y)\right|\lesssim \frac{1}{\mu(B(x,r^{1/m}))} \left(1+\frac{d(x,y)}{r^{1/m}}\right)^{-n-2N-\epsilon}$$
where $n$ is the homogeneous dimension of the space and $N$ and ``other dimension parameter'' due to the homogeneous type ($N\geq0$ could be equal to $0$).
 \item $({\mathcal A}_r)_{r>0}$ is a semigroup : for all $t,s>0$ ${\mathcal A}_s{\mathcal A}_t={\mathcal A}_{s+t}$.
\end{itemize}
To such a collection, we build the following operator : for $Q$ a ball
$$ B_Q^*(f) = f - {\mathcal A}_{r_Q^m}(f).$$
Let us check that our assumptions hold with $q_1=\infty$ and $q_0=1$. \\
Considering a ball $Q$ it comes that
\begin{align*}
 \left\|A_Q^*(f)\right\|_{L^\infty(Q)} & \lesssim \frac{1}{\mu(Q)} \sup_{x\in Q} \sum_{j\geq 0} \int_{C_j(Q)} \left(1+\frac{d(x,y)}{r_Q}\right)^{-n-2N-\epsilon} |f(y)| d\mu(y) \\
 & \lesssim \frac{1}{\mu(Q)} \sum_{j\geq 0} 2^{-(n+2N+\epsilon)j} \int_{C_j(Q)} |f(y)| d\mu(y) \\
 & \lesssim  \sum_{j\geq 0} \gamma_j  \left(\frac{1}{\mu(2^jQ)}\int_{C_j(Q)} |f| \right)
\end{align*}
with
$$\gamma_j \lesssim 2^{-(n+2N+\epsilon)j} \frac{\mu(2^jQ)}{\mu(Q)} \lesssim 2^{-(2N+\epsilon)j} \lesssim 2^{-\epsilon j}.$$
So Assumption (\ref{assum}) is satisfied.

\gb Then it remains us to check Assumption (\ref{ass:cont}). Indeed it corresponds to a local version of Proposition 2.6 of \cite{DY} for $K=2$.
Let us consider a ball $Q$ and another one $R\subset 2Q$, we have to estimate
$$ \left\| B_R^*(A_Q^*(f)) \right\|_{L^\infty(R)}.$$
With denoting by $r_R$ and $r_Q$ the radii and using the commutativity and the semigroup property of the operators, we have
$$ B_R^* A_Q^*(f) = {\mathcal A}_{r_Q}(f)- {\mathcal A}_{r_R+r_Q}(f).$$
Then as $r_R\leq 2r_Q$, Proposition 2.6 of \cite{DY} proves that this quantity is bounded by what we called $\|f\|_{BMO_1}$.
Consequently we deduce that Assumption (\ref{ass:cont}) is satisfied too.

\gb
{\bf Conclusion :} In the framework of \cite{DY}, our Assumptions (\ref{assum}) and (\ref{ass:cont2}) are satisfied. We can also apply Theorem \ref{thm} and obtain John-Nirenberg inequalities for the $BMO_q$ spaces with $q\in(1,\infty)$ and we recover the results of section 3 in \cite{DY}.

\subsection{Case of a semigroup associated to a second order divergence form operator}

The aim of this subsection is to compare the results of \cite{HM} about John-Nirenberg inequalities for BMO spaces associated to a divergence form operator to our ones.

\mb Let us first recall the framework of \cite{HM}. Consider the Euclidean space $X=\R^n$ and $A$ be an $n\times n$ matrix-valued function satisfying the ellipticity condition~: there exist two constants $\Lambda\geq \lambda>0$ such that
$$ \forall \xi,\zeta\in \C^n, \qquad \lambda|\xi|^2 \leq Re \left( A\xi \cdot \overline{\xi} \right) \quad  \textrm{and} \quad |A\xi \cdot \overline{\zeta} | \leq \Lambda|\xi||\zeta|.$$ We define the second order divergence form operator
$$ L(f):= -\textrm{div} (A \nabla f).$$
Semigroup associated to such operators satisfies to ``Gaffney estimates''~:

\begin{prop}[Lemma 2.5 \cite{HM} and \cite{A}] \label{prop:gaffney} There exist exponents $1\leq p_L<2<\widetilde{p_L}<\infty$ such that for every $p$ and $q$ with $p_L<p\leq q<\widetilde{p_L}$ the semigroup $(e^{-tL})_{t>0}$ satisfies to $L^p-L^q$ off-diagonal estimates, i.e. for arbitrary closed sets $E,F\subset \R^n$~:
 $$ \|e^{-tL}f\|_{L^q(F)} \lesssim t^{\frac{n}{2} \left(\frac{1}{q}-\frac{1}{p}\right)} e^{-\frac{d(E,F)^2}{t}} \|f\|_{L^p(E)}$$
for every $t>0$ and every function $f\in L^p(E)$.
\end{prop}

\mb In \cite{HM} S. Hofmann and S. Mayboroda define for $p\in(p_L,\tilde{p_L})$ a Hardy space $H^1_{L,p}$ associated to
this operator and give several charaterizations.
For $f\in L^1$ we have the equivalence of the following norms~:
\begin{align*}
 \|f\|_{H^1_{L,p}} & := \|f\|_{L^1}+ \left\| \left(\iint_{\genfrac{}{}{0pt}{}{t>0,\ y\in\R^n}{|x-y|\leq t}} \left|t^2L e^{-t^2L} f(y) \right|^2 \frac{dtdy}{t^{n+1}} \right)^{1/2} \right\|_{L^1} \label{maximaldef} \\
  & \simeq \|f\|_{L^1} + \left\| \sup_{\genfrac{}{}{0pt}{}{t>0,\ y\in\R^n}{|x-y|\leq t}} \left(\frac{1}{t^n} \int_{B(y,t)} \left|e^{-t^2L} f(z) \right|^2 dz \right)^{1/2} \right\|_{L^1}. \nonumber
  \end{align*}
In addition, they prove a molecular decomposition with the following definition : let $\epsilon>0$ and $M>n/4$ be fixed, a function $m\in L^2$ is a $(p,\epsilon,M)$-molecule if there exists a ball $Q\subset \R^n$ such that~:
\be{mol}
 \forall i\geq 0, \forall k\in\{0,...,M\}, \qquad \left\| \left(r_Q^{-2}L^{-1}\right)^{k} m \right\|_{L^p(C_i(Q))} \leq 2^{-i\epsilon} |2^{i+1}Q|^{-1+\frac{1}{p}}. \ee

\mb Moreover they prove that these spaces $H^1_{L,p}$ do not depend on $p\in (p_L,\widetilde{p_L})$ and they identify the dual spaces as a BMO space. For a ball $Q$, they consider operator $B_Q$ given by the radius of the ball by
$$ B_Q(f) = \left(I-e^{-r_Q^2L} \right)^{M}(f),$$
with a large enough integer $M>n/4$.

\gb
Let us check that our assumptions are satisfied in this case.\\
The operator $A_Q$ is given by the semigroup as follows
$$A_Q^*(f) := f-\left(I-e^{-r_Q^2L^*} \right)^{M}(f).$$
By expanding the power $M$, it comes that $A_Q^*$ is a finite combinaison of semigroups~:
$$A_Q^*(f) = \sum_{k=0}^M \left(\begin{array}{c}
                           M \\
                           k
                          \end{array} \right)  (-1)^{k} e^{-kr_Q^2L^*}(f).$$
So Gaffney estimates (see Proposition \ref{prop:gaffney}) give us some coefficients $\gamma_j$ (depending on $M$) such that Assumption (\ref{assum}) holds for $q_1=p_L'-\epsilon$ and $q_0=\widetilde{p_L}'+\epsilon$ with $\epsilon >0$ as small as we want.
It remains to check Assumption (\ref{ass:cont}), which is the goal of the following proposition (in fact we prove the stronger assumption (\ref{ass:cont2})).

\begin{prop} \label{prop:hm} In this framework, for $t>0$ the semigroup $e^{-tL^*}$ acts continuously in BMO spaces. Let write for a (as small as we want) parameter $\tau>0$, $q_1=p_L'-\tau$ and $q_0=\widetilde{p_L}'+\tau$.
\begin{itemize}
\item There is a constant $c$ such that for all $t>0$ and all exponent $q\in[q_0,q_1]$
$$ \left\| e^{-tL^*} f\right\|_{BMO_{q}} \leq c\|f\|_{BMO_{q}}$$
\item Assumption (\ref{ass:cont}) holds:
$$ \sup_{\genfrac{}{}{0pt}{}{R \textrm{ ball}}{r_R^2\leq 4t }} \left(\frac{1}{|R|}\int_{R} \left| B_R^*( e^{-tL^*} ) \right|^{q_1}  \right)^{1/q_1 } \leq c \|f\|_{BMO_{q_0}}. $$
\end{itemize}
\end{prop}

\dem Indeed we can prove a more precise result using duality. Thanks to Theorem 1.3 of \cite{HM}, the desired result is equivalent to the following~:
there is a constant $c$ such that for all $t>0$ and $p\in[q_1',q_0']$
$$ \left\| e^{-tL} f\right\|_{H^1_{L,p}} \leq c\|f\|_{H^1_{L,p}}$$
and
$$ \left\| e^{-tL} f\right\|_{H^1_{L,q_0'}} \leq c\|f\|_{H^1_{L,q_1',t}},$$
where $H^1_{L,q_1',t}$ is the Hardy space built on $(q_1',\epsilon)$ molecules associated to balls of radius lower than $2\sqrt{t}$.
This will be achieved by invoking the following lemma and the fact that the space $H^1_{L,q}$ does not depend on the integer $M$.
\findem

\begin{lem} \label{lem:hm} If $p\in[q_1',q_0']$ and $f$ is a $(p,\epsilon,2M)$-molecule adapted to a ball $Q$, then for every $t\geq 0$ $e^{-tL}(f)$ is a $(p,\epsilon, M)$ molecule adapted to the same cube. If $f$ is a $(q_1',\epsilon,2M)$-molecule adapted to a ball $Q$ (with $r_Q\lesssim \sqrt{t}$), then $e^{-tL}(f)$ is a $(q_0',\epsilon, M)$ molecule adapted to the same cube.
\end{lem}

\dem Let us first check the first claim. We need to show that $e^{-tL}(f)$ verifies (\ref{mol}) for $p$, up to some multiplicative constant that is uniform in $t$. We fix now the indices $i$ and $k$ and consider two cases. \\
\noindent {\bf Case 1:} $t\leq 2^i r_Q^2$. \\
If $i\leq 3$, we then have
\begin{align*}
 \left\|(r_Q^2L)^{-k} e^{-tL^2} (f) \right\|_{L^p(C_i(Q))} & = \left\|e^{-tL^2} (r_Q^2L)^{-k}  (f) \right\|_{L^p(C_i(Q))} \\
 & \lesssim \left\|(r_Q^2L)^{-k}  (f) \right\|_{L^p} \lesssim |Q|^{-1+\frac{1}{p}},
\end{align*}
as desired, where we have used the $L^{p}$-boundedness of the semigroup and the normalization of the molecule $f$ (\ref{mol}). \\
Suppose now $i>3$. We split
$$ (r_Q^2 L)^{-k}(f) = g_1+g_2$$
with
$$ g_1:=(r_Q^2 L)^{-k}(f) {\bf 1}_{2^{i-2}Q} \qquad g_2:=(r_Q^2 L)^{-k}(f) {\bf 1}_{(2^{i-2}Q)^c}.$$
We then have
\begin{align*}
 \left\|(r_Q^2L)^{-k} e^{-tL^2} (f) \right\|_{L^p(C_i(Q))} & = \left\|e^{-tL^2} (r_Q^2L)^{-k}  (f) \right\|_{L^p(C_i(Q))} \\
 & \leq \left\|e^{-tL^2} g_1\right\|_{L^p(C_i(Q))} + \left\|e^{-tL^2} g_2\right\|_{L^p(C_i(Q))} \\
 & \lesssim e^{-\frac{4^i r_Q^2}{t}} |Q|^{-1+1/p} + \sum_{j=i-1}^\infty 2^{-j(n/2+\epsilon)} |Q|^{-1+1/p}
\end{align*}
where in the last step we have used Gaffney estimate (Proposition \ref{prop:gaffney}), $L^{p}$-boundedness of the semigroup and (\ref{mol}). The desired bound follows in the present case.

\mb {\bf Case 2:} $t\geq 2^i r_Q^2$. \\
In this case, we have
\begin{align*}
 \left\|(r_Q^2L)^{-k} e^{-tL^2} (f) \right\|_{L^p(C_i(Q))} & = \left\|(r_Q^2L)^M e^{-tL^2} (r_Q^2L)^{-k-M}  (f) \right\|_{L^p(C_i(Q))} \\
 & = \left(\frac{r_Q^2}{t}\right)^M\left\|(tL)^M e^{-tL^2} (r_Q^2L)^{-k-M}  (f) \right\|_{L^p(C_i(Q))} \\
 & \lesssim 2^{-iM} |Q|^{-1+1/p}
\end{align*}
where in the last line we have used $L^{p}$-boundedness of $(tL)^M e^{-tL^2}$, along the fact that (\ref{mol}) holds with $k+M\leq 2M$ instead of $k$, since $f$ is a $(p,\epsilon,2M)$-molecule. Since we can choose $M>n/{p}$, the desired bound follows.

\mb
\noindent {\bf Second claim}. \\
It remains us to check the second claim, stated in the lemma. The proof is the same one as previously with the following modification. We have now to use the off-diagonal decay $L^{q_1'}-L^{q_0'}$ of the semigroup and the global boundedness (instead of the $L^{p}-L^{p}$ used before). Since Proposition \ref{prop:gaffney}, this operation makes appear an extra factor:
$$ \left(\frac{t^{n/2}}{|Q|}\right)^{\frac{1}{q_1'} - \frac{1}{q_0'}}.$$
The power $\frac{1}{q_1'} - \frac{1}{q_0'}$ is negative, so $\frac{t^{n/2}}{|Q|}$ should be bounded below in order than this new coefficient be bounded. That is why we require $r_Q\lesssim \sqrt{t}$.
\findem

\begin{rem} To prove that the semigroup operator $e^{-tL}$ continuously acts on the Hardy space $H^1_{L,p}$, we refer the reader to the work \cite{DP} of J.~Dziuba\'nski and M. Preisner. It is obvious that the function $x\rightarrow e^{-tx}$ satisfies to the above assumption of \cite{DP} and so the associated multiplier $e^{-tL}$ is bounded on the Hardy space (or at least on the molecules).
\end{rem}

\gb
{\bf Conclusion :} In the framework of \cite{HM}, our Assumptions (\ref{assum}) and (\ref{ass:cont2}) are satisfied. We can also apply Theorem \ref{thm} and obtain John-Nirenberg inequalities for the $BMO_q$ spaces with $q\in(\widetilde{p_L}',p_L')$. The precise inequality seems to be new, however we emphasize that the authors in \cite{HM} have already obtained an implicit John-Nirenberg inequalities in order to identify their BMO spaces making vary the exponent $q\in(\widetilde{p_L}',p_L')$ (see Section 10 in \cite{HM}).

\section{Application to Hardy spaces} \label{sec:appli}

We devote this section to an application of John-Nirenberg inequalities in the theory of Hardy spaces. We refer the reader to Subsection \ref{subsec:def} for definitions of {\it atoms} and {\it Hardy spaces}. We only deal with the atomic Hardy spaces for simplicity but a molecular version of the following results can be obtained too.

\mb First let us give a ``Hardy spaces''-point of view of our main Assumption (\ref{ass:cont2}).
\begin{rem}
Assumption (\ref{ass:cont2}) is equivalent to a $H^1_{p_1,ato}-H^1_{p_0,ato}$ boundedness of operators $A_Q$.
\end{rem}

\mb Now we assume that ${\mathbb B}$ satisfies some $L^p-L^p$ decay estimates : for $p\in[p_1,p_0]$
for $M''$ a large enough exponent, there exists a constant $C$ such
that \be{decay} \forall k\geq 0,\ \forall f\in L^p,\
\textrm{supp}(f) \subset  Q \qquad \left\| B_Q(f)
\right\|_{L^p(C_k(Q))} \leq C 2^{-M k} \|f\|_{L^p(Q)}. \ee

\noindent Using (\ref{decay}) we get the following properties about Hardy
spaces and BMO spaces.

\begin{prop} \label{inclusion} The Hardy space $H^1_{p,ato}$ is included in $L^1$ .
 \end{prop}

\dem Since its atomic decoposition, we only have to control the $L^1$-norm of each atom $m=B_Q(f_Q) \in H^1_{p,ato}$, by a uniform bound. \\
 By (\ref{decay}), the estimates for $f_Q$, the doubling property of $\mu$ and the fact that $M$ is large enough ($M >\delta/p'$ works), we have
\begin{align*}
\left\| B_Q(f_Q) \right\|_{L^1}   \leq  \sum_{k\geq 0} \left\|
B_Q(f_Q) \right\|_{L^1(C_k(Q))}
 \lesssim  \sum_{k\geq 0} \mu(Q)^{1/{p'}} 2^{-M k} \mu(Q)^{1/p -1 } 2^{{k\delta/p'}}  \lesssim 1.\\
\end{align*}
 So we obtain that each $p$-atom is bounded in $L^1$, which permits to complete the proof. \findem

\begin{cor}\label{Banach} For $p\in[p_1,p_0]$, the space $H^1_{p,ato}$ is a Banach space.
\end{cor}

\dem  The proof is already written in \cite{BZ}. We reproduce it here for an easy reference. \\
We only verify the completeness:
$H^1_{p,ato}$ is a Banach space if for all sequences
$(h_i)_{i\in\N}$ of $H^1_{p,ato}$ satisfying
$$\sum_{i\geq 0} \left\| h_i \right\|_{H^1_{p,ato}} <\infty, $$
the series $\sum h_i$ converges in the Hardy space
$H^1_{p,ato}$. \\
For such sequence in $H^1_{p,ato}$, we say that $\sum_{i} h_i \in L^1,$
 because each atom decomposition is absolutely
convergent in $L^1$-sense (since the previous proposition). If we denote $ f = \sum_{i} h_i \in L^1,$
then using the condition that $\sum_{i\geq 0} \left\| h_i
\right\|_{H^1_{p,ato}} <\infty, $ we have
$$ \| f - \sum_{i=0}^n h_i \|_{H^1_{p,ato}} \leq \sum_{i=n+1}^\infty \left\| h_i \right\|_{H^1_{p,ato}} \xrightarrow[n\rightarrow \infty]{} 0.$$ \findem

\mb We now come to our main result of this section: the Hardy space $H^1_{p,ato}$ does not depend
on $p\in(p_1,p_0]$.

\begin{thm} Under the above assumptions of Theorem \ref{thm} and (\ref{decay}), the Hardy space $H^1_{p,ato}$ does not depend on the exponent $p\in (p_1,p_0]$.
\end{thm}

\dem The proof is based on duality and the property that the BMO spaces are not depending on the exponent, see Corollary \ref{corco}. \\
We recall duality results and we refer to Section 8 of \cite{BZ} for more details. Fix an exponent $p\in (p_1,p_0]$. We cannot have a precise characterization of the dual space of our atomic Hardy space. However, we have the following results. Since the operator $B_Q$ are acting on a function (with a bounded support) to the Hardy space, we know that we can extend by duality $B_Q^*$ from $(H^1_{p,ato})^*$ to $L^{p'}_{loc}$. Then, we claim that for $\phi\in (H^1_{p,ato})^*$
\be{eq:duality}
 \left\| \phi \right\|_{BMO_{p'}} \simeq  \left\| \phi \right\|_{(H^1_{p,ato})^*}. \ee

\noindent {\bf First step} : Proof of (\ref{eq:duality}). \\
For each $ \phi \in BMO_{p'}\cap (H^1_{p,ato})^*,$ for each atom $ m \in
\left(H^1_{p,ato}\right), $ where $m = B_Q(f_Q), $ by H\"older
inequality, we have
\begin{equation}\label{eq100}
| \langle \phi, m  \rangle | = \left| \int_ Q B^*_Q (\phi) f_Q d\mu \right| \leq
\left(\int_ Q |B^*_Q (\phi) |^{p'}d\mu \right)^{\frac{1}{ p'}}
\left( \int_ Q |f|^p d\mu \right)^{\frac{1}{ p}} \leq \|
\phi\|_{BMO_{p'}}.
\end{equation}
Therefore by atomic decomposition, we deduce the first inequality $\left\| \phi \right\|_{(H^1_{p,ato})^*} \leq \left\| \phi \right\|_{BMO_{p'}} $. It remains us to check the reverse inequality. \\
For arbitrary  $f \in L^p(Q)$ satisfying $ \| f\|_{L^p(Q)} = 1$, we set $  g_Q := \mu (Q
)^{\frac{1}{p} -1}f, $ then $m=B_Q(g)$ is a $p$-atom.
Therefore
$$ \left| \int_ Q B^*_Q (\phi) f d\mu \right| = \mu (Q
)^{1-\frac{1}{p} } | \langle \phi, B_Q (g)
\rangle | \leq  \mu (Q
)^{\frac{1}{p'}} \|\phi \|_{(H^1_{p,ato})*}.
$$
This holds for every $L^p(Q)$-normalized function $f$. By duality, we deduce the reverse inequality, which concludes to (\ref{eq:duality}).

\mb {\bf Second step} : End of the proof. \\
Let chose two exponents $p,r$ in the above range. By symmetry, it is just sufficient to prove that
\be{eq:a} \|f\|_{H^1_{p,ato}}  \lesssim \|f\|_{H^1_{r,ato}}\ee
for every function $f\in H^1_{p,ato} \cap H^1_{r,ato}$ (since it is easy to check that $H^1_{p,ato} \cap H^1_{r,ato}$ is dense into both Hardy spaces).
So let us fix such a function $f$. The Hardy spaces beeing Banach spaces (Corollary \ref{Banach}),  Hahn-Banach Theorem implies that there is $\phi\in (H^1_{p,ato})^*$ normalized such that
$$  \|f\|_{H^1_{p,ato}} = \langle \phi, f \rangle.$$
We know that for every ball $Q$, $B_Q^*(\phi)$ belongs to $L^{p'}(Q)$ and satisfies
$$ \left(\frac{1}{\mu(Q)} \int_ Q |B^*_Q (\phi) |^{p'}d\mu \right)^{\frac{1}{p'}}\leq 1.$$
We can apply John-Nirenberg inequality: Theorem \ref{thm} \footnote{In Theorem \ref{thm}, we assume that $\phi$ is a function, which may not be the case here. However, we let the reader to check that the proof relies only  the fact that $B_Q^*(\phi)$ is measurable on $Q$.}.
We also obtain that for all ball $Q$, $B_Q^*(\phi)$ belongs to $L^{r'}(Q)$ and satisfies
$$ \left(\frac{1}{\mu(Q)} \int_ Q |B^*_Q (\phi) |^{r'}d\mu \right)^{\frac{1}{r'}}\lesssim 1.$$
Then as $f\in H^1_{p,ato} \cap H^1_{r,ato} \subset H^1_{r,ato}$, it follows by the previous reasonning (step 1) that
$$  \|f\|_{H^1_{p,ato}} = \langle \phi, f \rangle \lesssim \sum_i |\lambda_i  \langle \phi, m_i \rangle| \lesssim \|f\|_{H^1_{r,ato}},$$
where we have used an ``extremizing" decomposition $f= \sum_i \lambda_i m_i$ with $r$-atoms.
We also conclude to (\ref{eq:a}). \findem

% \subsection{Generalized Fefferman-Stein inequalities}
%
% We move the reader to \cite{BZ} and \cite{B} for the use of generalized Fefferman-Stein inequalities (in order to obtain interpolation results between Hardy spaces and Lebesgue spaces). In \cite{FS}, C.~Fefferman and E.M. Stein have introduced the maximal operator (based on the oscillation) and proved a maximal inequality, called {\it Fefferman-Stein inequality}. In \cite{BB}, the first author have extended these results in considering Sobolev versions.
%
% \gb In the framework of our abstract oscillation operators, (see Subsection \ref{subsec:def}), we define for an exponent $\sigma\in (1,\infty)$ the maximal sharp function
% $${\mathcal M}_\sigma^\sharp(f)(x) := \sup_{x\ni Q} \left(\frac{1}{\mu(Q)}\int_Q |B_Q^*(f)|^\sigma d\mu \right)^{1/\sigma}.$$
% Then, we ask the following question: for which exponents $q,\sigma\in(1,\infty)$, do we have the following equivalence
% \be{eq:FS} \|f\|_{L^q(X)} \simeq \|{\mathcam M}_\sigma(f) \|_{L^q}. \ee
%
% \begin{prop} If the space $X$ is of infinite measure then under the assumptions of Theorem \ref{thm}, (\ref{eq:FS}) holds for every exponents $q,\sigma \in(q_0,q_1)$.
% \end{prop}
%
% \dem  We do not repeat the proof, as such inequalities are a direct consequence of Assumption (\ref{assum}), due to Remark \ref{rem:int}. We just precise the new arguments. Following

\gb {\bf Acknowledgment:} We are indebted to Pr. Steve Hofmann for valuable advices and his help about Lemma \ref{lem:hm}.

\end{document}